\documentclass[leqno,12pt]{article}
\usepackage{amssymb,amsmath,graphics}
\title{ The equivariant Todd genus of a complete toric variety, with Danilov condition}
\author{Nicole Berline and Mich{\`e}le Vergne}
\date{October 2006}
\begin{document}
\maketitle
\newtheorem{theorem}{Theorem}
\newtheorem{proposition}[theorem]{Proposition}
\newtheorem{lemma}[theorem]{Lemma}
\newtheorem{definition}[theorem]{Definition}
\newtheorem{corollary}[theorem]{Corollary}
\newtheorem{remark}[theorem]{Remark}
\newtheorem{example}[theorem]{Example}
\newenvironment{proof}{{\bf Proof. }}{\par}
\newcommand{\C}{{\mathbb C}}
\newcommand{\R}{{\mathbb R}}
\newcommand{\Z}{{\mathbb Z}}
\newcommand{\N}{{\mathbb N}}
\newcommand{\CE}{{\cal E}}
\newcommand{\CF}{{\cal F}}
\newcommand{\la}{{\langle}}
\newcommand{\ra}{{\rangle}}
\newcommand{\mult}{\operatorname{mult}}
\newcommand{\divi}{\operatorname{div}}
\newcommand{\Supp}{\operatorname{Supp}}
\newcommand{\lin}{\operatorname{lin}}
\newcommand{\Todd}{\operatorname{Todd}}
\renewcommand{\a}{{\mathfrak{a}}}
\renewcommand{\c}{{\mathfrak{c}}}
\renewcommand{\d}{{\mathfrak{d}}}
\newcommand{\f}{{\mathfrak{f}}}
\newcommand{\g}{{\mathfrak{g}}}
\newcommand{\p}{{\mathfrak{p}}}
\renewcommand{\t}{{\mathfrak{t}}}
\newcommand{\lattice}{\Lambda}

{{\em Dedicated to Ernest Vinberg.}}

\begin{abstract}
 We write the \emph{equivariant} Todd class of a general
complete toric variety as an explicit  combination of the orbit
closures, the coefficients being  {\em analytic} functions on the
Lie algebra of the torus which satisfy Danilov's requirement.
\end{abstract}

\section{Introduction}
In \cite{EML}, to any rational affine cone $\a$ in a rational
vector space with a rational scalar product, we associated in a
canonical way a holomorphic function with rational Taylor
coefficients $\mu(\a)$. These functions are defined recursively by
an elementary combinatorial construction and $\mu(\a)$ can be
computed in polynomial time at any order, when $\dim \a$ is fixed.

In the present article, these  $\mu$-functions are used to write the
\emph{equivariant} Todd class of a general complete toric variety as
an explicit  combination of the orbit closures, the coefficients
being (germs of) {\em analytic} functions on $\g$ which satisfy
Danilov's requirement \cite{danilov}. Let $G\simeq(\C^\ast)^d$ be a
complex torus. Let $\sigma$ be a rational cone in the lattice of one
parameter subgroups of $G$, and let $\check{\sigma}$ be the
quotiented dual cone (see Definition \ref{mucech}). Then for any fan
$\Sigma$ in the lattice, the equivariant Todd class
$\Todd(X_\Sigma)$ of the corresponding $G$-toric variety $X_\Sigma$
is written as  the combination
\begin{equation}\label{localtodd}
\Todd(X_\Sigma)=\sum_{\sigma \in \Sigma}\mu(\check{\sigma}) [\bar{O}_\sigma]
\end{equation}
where $[\bar{O}_\sigma]$
is the equivariant homology class of the $G$-orbit closure $\bar{O}_\sigma$.
Notice the special features of this formula:
 all orbits closures $\bar{O}_\sigma$ appear and, for each $\sigma$,
the term $\mu(\check{\sigma})$ is a germ of  analytic function (with
rational Taylor coefficients) on $\g$ \emph{which depends only on
the cone $\sigma$, not on the fan}
 (Danilov condition).

Let us write (\ref{localtodd}) in the case of the projective line
$\mathbb{P}_1(\C)$. There are three orbit closures,
$\mathbb{P}_1(\C)$ itself  corresponding to the cone $\{0\}$, and
the north and south poles $p_+$ , $p_-$ corresponding respectively
to the  cones $\R_+$ and $-\R_+$.

Define
$$
B(\xi)= \frac{1}{e^\xi-1}- \frac{1}{\xi}= \sum_{n\geq 1}\frac{b(n)}{n!}\xi^{n-1},
$$
where $b(n)$ are the Bernoulli numbers.
We obtain the following expression for the equivariant Todd class of $\mathbb{P}_1(\C)$:
\begin{equation}\label{P1}
\Todd(\mathbb{P}_1(\C))(\xi)= [\mathbb{P}_1(\C)] - B(\xi)[p_+]-
B(-\xi)[p_-].
\end{equation}

\bigskip

The construction of this paper extends to equivariant homology a
similar result obtained by  Pommersheim-Thomas \cite{pomm} for the
ordinary Todd class. The construction of \cite{pomm} requires a
choice of a complement map. In particular, a scalar product provides
a natural complement map. Previously, Morelli \cite{morelli} gave
another type of construction for the ordinary Todd class, satisfying
Danilov condition, depending (rationally) of a choice of a flag in
$\g$. Other choices of complement maps, encompassing Morelli's
construction, are studied in \cite{pomm}.

In \cite{BriVer97}
 Brion and Vergne express the Todd class
 in \emph{equivariant cohomology}, in the case when the associated
 fan is simplicial, as a linear combination of products of the canonical
invariant divisors of the toric variety. These divisors generate the
equivariant cohomology as an algebra.

In \cite{brylinskizhang},  Brylinski and  Zhang give a  formula
for the equivariant Todd class of  a general complete toric
variety, as an element of a localized equivariant \emph{homology}
module. The Todd class is written as a linear combination of  the
homology classes defined by fixed points, the coefficients being
germs of {\em meromorphic} functions on $\g$. Their result can be
recovered from our formula (\ref{localtodd}) by using the
relations between the equivariant homology classes of orbit
closures $[\bar{O}_\sigma]$, (\cite{brylinskizhang}, Lemma 8.8).
For instance, in the case of the projective line
$\mathbb{P}_1(\C)$, one has the relation
$$
-\xi\; [\mathbb{P}_1(\C)]= [p_+]- [p_-].
$$
From this relation and (\ref{P1}), we recover the formula of
Brylinski-Zhang for  $\mathbb{P}_1(\C)$:
$$
\Todd(\mathbb{P}_1(\C))(\xi)=\frac{1}{1-e^{\xi}}\; [p_+]+
\frac{1}{1-e^{-\xi}}\; [p_-]
$$
In order to compute the Todd class, we start with the results of
\cite{BriVer97}. In the case of a simplicial variety,  we rewrite
the Todd class in cohomology using our  Euler-Maclaurin formula of
\cite{EML}, then we translate it in homology using equivariant
Poincar{\'e} duality for a simplicial toric variety. Finally, for
the general case, we use the push-forward method in equivariant
homology devised in \cite{brylinskizhang}.

\bigskip
When the fan $\Sigma$ is dual to an integral polytope $\p$,
Formula (\ref{localtodd}) is a rewriting, via Riemann-Roch, of the  local Euler-Maclaurin expansion
obtained in \cite{EML} for the sum $\sum_{x\in \p\cap\Z^d}h(x)$ of the values $h(x)$ at
integral points of $\p$, where
$h(x)$ is a polynomial function on $\R^d$.  This sum is written as a
sum of integrals over the various faces of $\p$,
\begin{equation}\label{intro-maclaurin}
\sum_{x\in \p\cap\Z^d}h(x)=\sum_{\f\in\CF(\p)}\int_\f
D(\p,\f)\cdot h
\end{equation}
where,  for  each face $\f$ of $\p$ , $D(\p,\f)$ is the differential
operator (of infinite order) with constant coefficients on $\R^d$
whose symbol is the function $\mu(\t(\p,\f))$ associated  to the
transverse cone of $\p$ along $\f$.
In dimension one, the partner of (\ref{P1}) is the classical Euler-Maclaurin summation formula
for a lattice interval $[a_1,a_2]$
\begin{eqnarray*}
\sum_{a_1}^{ a_2}h(x) = \int_{a_1}^{a_2}h(t)dt -  \sum_{n\geq
1}\frac{b(n)}{n!}h^{(n-1)}(a_1) +\sum_{n\geq 1}(-1)^{n}
\frac{b(n)}{n!} h^{(n-1)}(a_2).
\end{eqnarray*}

\medskip
Let us note that the combinatorial formula (\ref{intro-maclaurin})
of \cite{EML} holds for all rational polytopes, not only integral
ones.

Articles on toric varieties often stress their application to  lattice points in polytopes.
In this study, the relationship goes the other way round,
a situation not uncommon in the interaction between science and technology:
a progress in applications leads to a theoretical progress.

\section{Definitions and notations}

\subsection{Cones}
Let $V$ be a rational vector space, that is to say  a finite
dimensional real vector space with a lattice denoted by
$\lattice_V$ or simply $\lattice$. We denote by $V^* $ the dual
space of $V$.  We will denote elements of $V$ by latin letters
$x,y,v,\dots$ and elements of $V^*$ by greek letters $\xi,\dots$.
We denote the duality bracket by $\langle\xi,x\rangle$.

$V^* $ is equipped with the dual lattice $\lattice^* $
of $\lattice$ :
$$
\lattice^* = \{\xi\in V^*\; ; \,\,\langle\xi,x\rangle \in \Z \; \mbox{for
  all}\;
  x\in \lattice \}.
$$
If $S$ is a subset of $V$, we denote by $\lin(S)$ the vector
subspace of $V$ generated by  $S$ and by  $S^{\perp}$ the subspace
of $V^*$ orthogonal to $S$:
$$
S^{\perp}= \{\xi\in V^*\; ;\,\,\langle\xi,x\rangle =0 \;\mbox{for
all}\;
  x\in  S\}.
$$
A subspace $W$ of $V$ is called rational if $W\cap \lattice $ is a
lattice in $W$.  If $W$ is a  rational subspace, the image of
$\lattice$ in $V/W$ is a lattice in $V/W$, so that $V/W$ is a
rational vector space.

The space  $V$, with lattice $\lattice$, has a canonical Lebesgue
measure, for which  $V/\lattice$ has measure $1$.

The set of non  negative real numbers is denoted by $\R_+$. A convex
rational cone  $\c$ in a  rational space is  a closed convex cone
$\sum_{i=1}^k\R_+ v_i$ which is generated
 by a finite number of lattice vectors.
In this article, we simply say  cone instead of convex rational
cone.

 A  cone $\c$ is called simplicial if it is  generated by linearly
independent primitive vectors $v_1,\dots, v_k$. The multiplicity
$\mult(\c)$ of the simplicial cone $\c$ is the cardinal of the set
$\Lambda\cap \lin(\c)/\sum_{i=1}^k\Z v_i$. It is equal to
$|\det(v_1,\dots, v_k)|$ where the determinant is computed with
respect to a basis of the lattice $\Lambda\cap \lin(\c)$.

The set of faces of a cone $\c$ is denoted by $\CF(\c)$. If $\c$ is
pointed, then its vertex $0$  is its unique face of dimension $0$,
while $\c$ is the unique face of maximal dimension $\dim \c$.

If $\f$ is a face of the cone $\c$, the  \emph{transverse cone of
$\c$ along $\f$} is the image of $\c$ in the quotient space
 $V/\lin(\f)$. We denote it by $\t(\c,\f)$.

 Let $\c$ be a cone in $V$. The dual cone $\c^*$ of
$\c$ is the set of $\xi \in V^*$ such that $\la \xi,x\ra \geq 0$
for any $x \in \c$.

\medskip

We shall make use of subdivisions of cones.
\begin{definition}
A subdivision of a cone $\c$ is a finite collection $\cal C$ of
cones in $\lin(\c)$ such that:

(a) The faces of any  cone  in $\cal C$ are in $\cal C$.

 (b) If $\d_1$ and $\d_2$ are two elements of $\cal C$, then the
 intersection $\d_1\cap \d_2$ is a face  of both  $\d_1$ and $\d_2$.

 (c) We have  $\c=\cup_{\d \in \cal C}\d.$

The subdivision is called simplicial if it consists of  simplicial
cones.
\end{definition}

\begin{example}\label{subdivisionR}
The basic example is the subdivision $\{\R_+,\R_-, \{0\}\}$ of the
one-dimensional cone $\R$ where $\R_-$ denotes the opposite cone
$-\R_+$.
\end{example}

\subsection{A holomorphic function  associated to a cone}\label{mainconstruction}

We now recall the construction of \cite{EML}. Let $V$   be  a
rational space with lattice $\lattice_V$. We fix a scalar product
$Q(x,y)$ on $V$.  We assume that $Q$ is rational, meaning that
$Q(x,y)$ is rational for $x,y  \in \lattice_V$. To any cone $\c$
in a quotient space $W$ of $V$, we associate in \cite{EML} a (germ
at $0$ of) holomorphic function on $W ^*$, denoted by
$\mu(\c)(\xi)$.

The functions $\mu(\c)(\xi)$ satisfy  the following \emph{local
Euler-Maclaurin formula}.
\begin{equation}\label{definitionmu}
S(\c)(\xi)= \sum_{\f\in \CF(\c) }\mu(\t(\c,\f))(\xi)I(\f)(\xi),
\end{equation}
for $\xi\in W^*$ small.

 This formula requires some explanations
(see \cite{EML} for details).

The  function $S(\c)(\xi)$ is defined as a meromorphic function on
$W^*$ as follows. First, if the cone $\c$ is pointed, the series
$\sum_{x\in \c \,\cap \lattice_W} e^{\la \xi,x\ra}$ converges only
for those $\xi \in W^*$ such that $\la \xi,x\ra <0$ for all $x\in \c
\backslash \{0\}$. The sum can be extended to a
 meromorphic function on $W^* $. This is $S(\c)(\xi)$.
 If $\c$ is not pointed, the sum converges nowhere, and $S(\c)(\xi)\equiv 0$ by definition.

The integral $I(\f)(\xi):= \int_\f e^{\la \xi,x\ra}dm_\f(x)$  is
defined as a  rational function on $W^*$ in a similar way. The
measure $dm_\f(x)$ is the normalized Lebesgue measure on
$\lin(\f)$ defined by the lattice $\lattice_W \cap \lin(\f)$.

In Formula (\ref{definitionmu}), the transverse cone $\t(\c,\f)$ is
a cone in the quotient space $W/\lin(\f)$. The function
$\mu(\t(\c,\f))$ is a function on a neighborhood of $0$ in the dual
$(W/\lin(\f))^*\cong \lin(\f)^\perp\subset W^*$. We give a sense to
$\mu(\t(\c,\f))(\xi)$ for a small  $\xi\in W^*$ by extending this
function to a neighborhood of $0$ in the space $W^*$ itself by means
of orthogonal projection $W^*\to \lin(\f)^{\perp}$.

We have
$$\mu(\{0\})= 1   \mbox{ and } \mu(\R_+)(\xi)=\frac{1}{1-e^{\xi}}+\frac{1}{\xi}.$$
\begin{remark}
In \cite{EML} the function $\mu(\a) $ is defined for any rational
\emph{affine }cone $\a$. Here we will need only the case of a cone
with vertex $0$.
\end{remark}

\begin{remark} \label{after} The sum $S(\c)$ is easy to compute when $\c$ is a simplicial cone in
$V$ of dimension equal to $d=\dim V$. Let $v_1,v_2,\ldots,v_d$ be
the primitive vectors of the edges of $\c$. We denote by $
\Box(\c)=\sum_{i=1}^{d}[0,1[v_i$ the semi-open parallelepiped
generated by  the $v_i$'s .
 Then
\begin{equation}\label{Sparallelepiped} S(\c)(\xi)=
\left(\sum_{x\in \Box(\c)\cap \lattice_V} e^{\la
  \xi,x\ra} \right)\prod_{i=1}^d \frac{1}{1- e^{\la \xi,v_i\ra}}.
\end{equation}

In particular, the function  $\xi\mapsto \prod_{i=1}^d \la
\xi,v_i\ra S(\c)(\xi)$  is  holomorphic  near $0$ in $V^*$.
\end{remark}

Recall the relation between the rational function $I(\c)$ and the
rational functions $I(\f)$ where $\f$ varies over the set of
facets of $\c$.

Let $\c$ be a pointed cone in $V$ such that $\lin(\c)=V$. Let
$\c^*$ be its dual cone in $V^*$, and let $\CE(\c^*)$ be the set
of its edges. An edge $\tau$ of $\c^*$ defines the facet
$\f_\tau=\lin(\tau)^{\perp}\cap \c$ of $\c$. We choose primitive
generators $n_\tau$ on each edge $\tau$.

\begin{lemma}\label{relationsstokes}
For $\xi\in V^*$ and $v\in V$, we have
$$-\la v,\xi\ra I(\c)(\xi)=\sum_{\tau\in \CE(\c^*)}\la v,n_\tau\ra
I(f_\tau)(\xi).$$
\end{lemma}
\begin{proof}
 Let us denote by $\lambda$ the
Lebesgue volume form on $V$ defined by the lattice and an
orientation. We have, for $\c$ oriented accordingly,
$$
-\la v,\xi\ra I(\c)(\xi)=-\int_ {\c} e^{\langle\xi,x\rangle} \la
v,\xi\ra\lambda=-\int_\c d\alpha,
$$
where the differential form $\alpha$ on $V$  is given by
$$
\alpha = e^{\langle\xi,x\rangle}\imath_v \lambda.
$$
Then  Lemma \ref{relationsstokes}  is just  Stokes formula. $\Box$
\end{proof}

\subsection{Todd measure of a cone.}

Let $G\simeq(\C^\ast)^d$ be a complex torus. We denote its
character group by $\lattice$, the  real vector space spanned by
$\lattice$ by $\g^*$, the dual of $\g^*$ by $\g$ and the dual
lattice by $\lattice^*\subset \g$.

We denote by $S(\g^*)$ the symmetric algebra over $\g^*$, identified
with the algebra of polynomial functions on $\g$ and by
$\hat{S}(\g^*)$ its completion, the algebra of formal power series.
If $h$ is a holomorphic function on $\g$ defined near $0$, its
Taylor series is  an element of $\hat S(\g^*)$ also denoted by $h$.

In the following, we will  apply the construction of the preceding
paragraph \ref{mainconstruction} with $\g^*$ in the role of the
space $V$, therefore the dual $V^*$ will be the space $\g$.

Let $\sigma\subset V^*$ be a cone. The dual cone $\sigma^*$ of
$\sigma$ is the set of $x \in V$ such that $\la \xi,x\ra \geq 0$ for
any $\xi \in \sigma$.
 The vector subspace $\lin(\sigma)^\perp\subset
V$ of elements $x \in V$ such that $\la \xi,x\ra = 0$ for any $\xi
\in \sigma$ is contained in $\sigma^*$.
\begin{definition}\label{conecech}
Let $\sigma $ be a cone in $V^*$. We denote by $\check{\sigma}$ the
image of $\sigma^*$ in the quotient space $V/\lin(\sigma)^\perp$.
\end{definition}
The dual of $ \g^*/\lin(\sigma)^\perp$ is identified with the
subspace $\lin(\sigma)\subset \g$.
\begin{definition}\label{mucech}
 $\mu^*(\sigma)$ is  the (germ at $0$ of) holomorphic function on $\g$ defined as follows.
 For $\xi\in\lin(\sigma)$,
 $$
{\mu^*}(\sigma)(\xi):= \mu(\check{\sigma})(\xi),
 $$
 then  $\mu^*(\sigma)$ is extended to $\g$ by orthogonal projection.
\end{definition}
A crucial property of these functions is that the assignment
$\sigma \mapsto {\mu^*}(\sigma)(\xi)$ is a \emph{simple
valuation}:
\begin{theorem}\label{valuation}(\cite{EML}, Corollary 21)
Let $\sigma_0$ be a cone in $V^*$, and let $\cal C$ be a subdivision
of $\sigma_0$. Then
$$
\mu^*(\sigma_0)= \sum_{\sigma\in {\cal C}; \,\dim \sigma = \dim
\sigma_0} \mu^*(\sigma).
$$
\end{theorem}
Thus, extending a definition of  \cite{pomm},  we may call
${\mu^*}(\sigma)(\xi)$ the \emph{equivariant Todd measure of
$\sigma$}.

\subsection{The G-toric variety defined by $\Sigma$ }

Let $G\simeq(\C^\ast)^d$ be a complex torus.
 Let
$\Sigma$ be a $G$-fan. For $1\leq i\leq d$, we denote by
$\Sigma[i]$ the set of  cones in $\Sigma$ of dimension $i$.

 The G-toric variety defined by
$\Sigma$ (\cite{danilov}, \cite{fulton}) is denoted by $X_\Sigma$.
Recall that a cone $\sigma$  in $\Sigma$  defines an open
$G$-invariant affine subvariety $U_\sigma$ of $X_\Sigma$. In this
article,  we follow the conventions of \cite{BriVer97}. Thus, if
$r\in \sigma^*\cap \Lambda$, then $r$ is a regular function on
$U_\sigma$ of weight $r$ with respect to the action of $G$ on
$X_\Sigma$. This convention determines the action of $G$ on
$X_\Sigma$. Note that this convention is contrary to that of
\cite{brion97}.

The open set $U_\sigma$ contains a unique relatively closed
$G$-orbit $O_\sigma$ whose closure in $X_\Sigma$ is denoted by
$\bar{O}_\sigma$. If $\sigma\in \Sigma[d]$, then $U_\sigma$ contains
a  unique $G$-fixed point $p_\sigma$  and
 $O_\sigma=\{p_\sigma\}$.

\begin{example}\label{THE}
Let $G=\C^*$. Denote by $\ell$ the  canonical generator of $\Lambda$
($\ell(u)=u$, for $u\in \C^*$). The projective line $P_1(\C)$ is
associated to the fan $,\R_-, \{0\}, \R_+$ in $\g$.
 The fixed points are  $p_+$ and  $p_-$,  the action of $G$ on the
tangent space at $p_+$ being given by $-\ell$, and the action of $G$
on the tangent space at $p_-$ being given by $\ell$. Then
$$\bar{O}_{\{0\}}=P_1(\C),\hspace{1cm}\bar{O}_{\R_+}=\{p_+\},
\hspace{1cm}\bar{O}_{\R_-}=\{p_-\}.$$
\end{example}

Denote by $R_\Sigma$ the algebra of continuous \emph{piecewise
polynomial} functions  on $\cup_{\sigma\in \Sigma}\sigma$. An
element of $R_\Sigma$ is a function $f$  such that, for any cone
$\sigma$ in $\Sigma[d]$,  the restriction of $f$ to $\sigma$  is
equal to the restriction to $\sigma$ of an element of  $S(\g^*)$
which we denote by $f|_\sigma$.

$R_\Sigma$ is a module over $S(\g^*)$. We define $\hat R_\Sigma=\hat
S(\g^*)\otimes_{S(\g^*)}R_\Sigma$. If $\sigma\in \Sigma[d]$, the
restriction map $f\mapsto f|_\sigma$ extends to a map from $\hat
R_\Sigma$ to $\hat S(\g^*)$.

Let $H_*^G(X_\Sigma) $ denote the equivariant homology group  of
$X_\Sigma$. It is a module over $S(\g^*)$, with the following
relations (see \cite{brion97}).  Any $G$-invariant closed
subvariety $Y$ of $X_\Sigma$ defines an element of
$H_*^G(X_\Sigma)$, denoted by $[Y]$. Let $v$ be an element of
$\Lambda$, that is a character of $G$. If $f$ is a function on $Y$
of weight $v$, then $-v[Y]=[\divi(f)]$.

\begin{example}
Return to Example \ref{THE}.
  Then
$-\ell[P_1(\C)]=[p_+]-[p_-]$.

Indeed, $\ell$ is a rational function $f$ on $P_1(\C)$ of weight
$\ell$ and $\divi(f)=\{p_+\}-\{p_-\}$
\end{example}

\section{A formula for the equivariant Todd class of a complete  simplicial
toric variety.}

In this section we assume that  $\Sigma$ is a complete simplicial
fan. Then $X_\Sigma$ is complete and is locally the quotient of a
smooth variety by a finite group. One says that   $X_\Sigma$ is a
simplicial complete toric variety.

Let $H^*_G(X_\Sigma)$ be the equivariant cohomology algebra  of
$X_\Sigma$. The Todd class ${\rm Todd}_{H^*_G}(X_\Sigma)$ is  an
element of the completed
 algebra $\hat{H}_G^*(X_\Sigma)$, where

$$\hat{H}_G^*(X_\Sigma)= \hat{S}(\g^*)\otimes_{S(\g^*)}H_G^*(X_\Sigma).$$

Let us recall Brion's isomorphism between $H^*_G(X_\Sigma)$ and
$R_\Sigma$.

 Let $\alpha\in H^*_G(X_\Sigma)$.
For every  $\sigma\in \Sigma[d]$, the restriction of $\alpha$ to the
fixed point $p_\sigma$ is denoted is an element of $S(\g^*)$ which
we denote by by $\alpha(p_\sigma)$.
\begin{theorem}\cite{brion}
There exists a unique algebra  isomorphism $B$ of $H_G^*(X_\Sigma)$
with $R_\Sigma$ such that, for $\alpha\in H_G^*(X_\Sigma)$ and
$\sigma\in \Sigma[d]$, one has $B(\alpha)|_\sigma=\alpha(p_\sigma)$.
\end{theorem}

This isomorphism $B$ extends to an isomorphism from
$\hat{H}_G^*(X_\Sigma)$ to $\hat{R}_\Sigma$.
\bigskip

We recall the formula (\cite{BriVer97}, 4.1, Theorem) which
expresses the Todd class ${\rm Todd}_{H^*_G}(X_\Sigma)$ of
$X_\Sigma$, as an element of $\hat{R}_\Sigma$

Let $\sigma\in \Sigma[d]$. The dual cone $\sigma^*$ is a simplicial
cone of dimension $d$ in $\g^*$. Let $v_i$ be the primitive
generators of the edges of $\sigma^*$.
\begin{equation}\label{toddsimplicial}
B({\rm
Todd}_{H^*_G}(X_\Sigma))|_\sigma(\xi)=\frac{1}{\mult(\sigma^*)}\prod_{i=1}^d
\la -\xi,v_i\ra S(\sigma^*)(\xi)   \mbox{  for  }\xi\in\g.
\end{equation}

\begin{remark}
By Remark \ref{after}, the right hand side is indeed holomorphic
at $0$.
\end{remark}

\bigskip

Let $\sigma$ be any cone in  $\Sigma$. Following (\cite{BriVer97},
3.3), we consider the following element $\varphi_\sigma$ of
$R_\Sigma$.

Let $\tau\in \Sigma[1]$ be a cone of dimension $1$ and let
$\eta_\tau$ be the generator of $\tau\cap \Lambda^*$. We denote by
$\varphi_\tau\in R_\Sigma$ the continuous piecewise linear
function on $\g$ such that $\varphi_\tau(\eta_\tau)=1$ and
$\varphi_\tau(\eta_{\tau'})=0$ for all $\tau'\in \Sigma[1]$
different from $\tau$.

\begin{definition}
Let $\sigma$ be a cone in $\Sigma$ and $mult(\sigma)$ its
multiplicity. Let $\CE(\sigma)$ be the set of edges of $\sigma$.
Then we set:
$$\varphi_\sigma=\mult(\sigma)\prod_{\tau\in \CE(\sigma)}\varphi_\tau.
$$
\end{definition}

The function $\varphi_\sigma$ vanishes identically on all cones
which do not contain $\sigma$.

We compute the restriction  $\varphi_\sigma|_{\sigma_0}$ for  a cone
$\sigma_0$ of dimension $d$. If $\sigma_0$ does not contain
$\sigma$, then $\varphi_\sigma|_{\sigma_0}=0$. Otherwise, let
$\sigma_0\in \Sigma[d]$ containing $\sigma$, let $
\eta^{\sigma_0}_1,\dots,\eta^{\sigma_0}_d$ be primitive  generators
of the edges of $\sigma_0$ such that $
\eta^{\sigma_0}_1,\dots,\eta^{\sigma_0}_r$ generate $\sigma$, and
let $v^{\sigma_0}_1,\dots, v^{\sigma_0}_d $ be the primitive
generators  of $\sigma_0^*$, indexed in such a way that $\la
v^{\sigma_0}_i,\eta^{\sigma_0}_i \ra>0$. Let $\f$ be the cone in
$\g^*$ generated by $v^{\sigma_0}_{r+1},\dots,v^{\sigma_0}_d$. Thus
$\lin(\f)=\lin(\sigma)^{\perp}$.
\begin{lemma}\label{check}
$$\varphi_\sigma|_{\sigma_0 }=\frac{\mult(\f)}{\mult(\sigma_0^*)}
\prod_{j=1}^r  v^{\sigma_0}_j.
$$
\end{lemma}

\begin{proof}
Let $q_j=\la v^{\sigma_0}_j,\eta^{\sigma_0}_j \ra$, so that the
dual basis to the basis $\{\eta^{\sigma_0}_j\}_{j=1}^d $ consists
of the vectors $\{\frac{1}{q_j}v^{\sigma_0}_j\}_{j=1}^d$. The
definition of $\varphi_\sigma$ implies that
$$\varphi_\sigma|_{\sigma_0}=\mult(\sigma)\prod_{j=1}^r\frac{1}{q_j}v^{\sigma_0}_j.$$
Thus we need to prove:
$$\frac{\mult(\f)}{\mult(\sigma_0^*)}=\mult(\sigma)\prod_{j=1}^r\frac{1}{q_j}.$$

This equality is obtained immediately by computing the
multiplicities as absolute values of  determinants, using a basis of
$\Lambda$  in which the matrix of the vectors $\eta_1,\dots,\eta_d$
is upper triangular. $\Box$
\end{proof}

As $\Sigma$ is a complete simplicial fan, the Poincar{\'e}
isomorphism $P: H_*^G(X_\Sigma)\mapsto H_G^*(X_\Sigma)$ is an
isomorphism  of $S(\g^*)$-modules between the equivariant homology
and the equivariant cohomology.

\begin{proposition}\label{FormulaP}(\cite{BriVer97}, 3.3)
$$
B(P({\bar O}_\sigma))=(-1)^{\dim \sigma}\varphi_\sigma.
$$
\end{proposition}

The classes $\varphi_\sigma$, when $\sigma$ varies over all cones
in $\Sigma$,  generates $R_\Sigma$ as a module over $S(\g^*)$.
Thus the equivariant Todd class ${\rm Todd}_{H^*_G}(X_\Sigma)$ of
a complete simplicial toric variety can be expressed as a
 combination of the classes $\varphi_\sigma$ with coefficients in
 $\hat S(\g^*)$. Such an expression is highly not unique as there are
 relations between the elements $\varphi_\sigma$.
  From now on,
  we choose a rational scalar product on $\g^*$.
  Given this choice, we are able to give
 a canonical formula.

\begin{proposition}\label{expansionsimplicial}
Let $G$ be a torus with Lie algebra $\g$. Let $\Sigma$ be a
complete rational simplicial fan in $\g$ and let $X_\Sigma$ be the
corresponding $G$-toric variety. Choose a scalar product on $\g$.
Then the equivariant Todd class ${\rm Todd}_{H^*_G}(X_\Sigma)$ of
$X_\Sigma$ is given by the following  combination of the classes
$\varphi_\sigma$ with coefficients in $\hat S(\g^*)$:
 \begin{equation}\label{tt}B({\rm Todd}_{H_G^*}(X_\Sigma))=
\sum_{\sigma\in \Sigma}(-1)^{\dim
\sigma}\mu^*(\sigma)\varphi_\sigma.
\end{equation}
\end{proposition}

 \begin{example}
Return to Example \ref{THE}. The restriction of the Todd class to
$\R_+$ is $\frac{-\xi}{1-e^{\xi}}$. Compute the restriction to
$\R^+$ of the right hand side of Formula (\ref{tt}):
$$\mu^*(\{0\})\varphi_{\{0\}}-\mu^*(\R_+)\varphi_{\R_+}
-\mu^*(\R_-)\varphi_{\R_-}.$$
 We obtain
$$1-\xi(\frac{1}{1-e^{\xi}}+\frac{1}{\xi})=\frac{-\xi}{1-e^{\xi}}.$$

In the same way, we check the equality on $\R_-$.
\end{example}
\begin{proof}
We only need to check that both members of Equation (\ref{tt}) agree
on each  cone $\sigma_0\in \Sigma[d]$. Let
$\eta_1,\eta_2,\ldots,\eta_d$ be the primitive generators of the
edges of $\sigma_0$. Let $v_1, v_2,\ldots, v_d$ be the primitive
generators of the edges of $\sigma_0^*$, so that $\la
v_i,\eta_i\ra>0$.  As $\varphi_\sigma$ vanishes if $\sigma$ is not
contained in $\sigma_0$, the right hand side of Equation (\ref{tt})
restricted to $\sigma_0$ is a sum over the faces of $\sigma_0$.
  For any subset $J$ of
$\{1,\dots,d\}$, we denote by $\sigma_J\subset\sigma_0 $ the cone
generated by the vectors $\eta_j$ for $j\in J$. The cone $\{0\}$
corresponds to the empty set.

Let us write the local Euler-Maclaurin formula (\ref{definitionmu})
for the cone $\sigma_0^*$. We label the faces $\f_J$ of $\sigma_0^*$
by the subsets $J\subset\{1,\dots,d\}$, by setting
 $\f_J=\sum_{j\notin J}\R_+ v_j$. Thus, for each subset
$J$, we have
 $\lin(\f_J)=\lin(\sigma_J)^\perp.$ The transverse cone $\t(\sigma_0^*,\f_J)$ and the projected cone
$\check{\sigma}_J \subset \g^*/\lin(\sigma_J)^\perp$ coincide. The
integral $I(\f_J)(\xi)$ is given by
$$
I(\f_J)(\xi)= \mult(\f_J)\prod_{j\notin J}\frac{1}{\langle
-\xi,v_j\rangle}.
$$
We obtain
\begin{equation}\label{EMLsimplicial}
S(\sigma_0^*)(\xi)= \sum_J \mult(\f_J)\mu(\check{\sigma}_J
)(\xi)\prod_{j\notin J}\frac{1}{\langle -\xi,v_j\rangle}.
\end{equation}
Thus \begin{eqnarray*} \frac{1}{\mult(\sigma_0^*)}\prod_{1\leq
j\leq d}\langle -\xi,v_j\rangle S(\sigma_0^*)(\xi)&=& \sum_J
\frac{\mult(f_J)}{\mult(\sigma_0^*)}\mu(\check{\sigma}_J )(\xi)
\prod_{j\in J}{\langle -\xi,v_j\rangle}\\
&=&\sum_J
(-1)^{|J|}\mu(\check{\sigma}_J)(\xi)\frac{\mult(\f_J)}{\mult(\sigma_0^*)}
\prod_{j\in J}{\langle \xi,v_j\rangle}\\
&=&\sum_J
(-1)^{|J|}\mu(\check{\sigma}_J)(\xi)\varphi_{\sigma_J}(\xi).
\end{eqnarray*}
The last equality follows from Lemma \ref{check}. Hence, recalling
Equation \ref{toddsimplicial}, we have proven Proposition
\ref{expansionsimplicial}. $\Box$
\end{proof}

\section{A formula for the Todd class of a complete toric variety}

In this section we  prove the formula (\ref{localtodd}) announced in
the introduction, in the  case of any complete toric variety, not
necessary simplicial.

Let $\Sigma$ be a complete fan in $\g$. We do not assume that
$\Sigma$ is  simplicial. The equivariant Todd class ${\rm
Todd}_{H^G_*}(X_\Sigma)$ of $X_\Sigma$ is then defined as an
element of the completed equivariant \emph{homology} module
$\hat{H}_*^G(X_\Sigma)$,

 $$\hat{H}_*^G(X_\Sigma)= \hat{S}(\g^*)\otimes_{S(\g^*)}H_*^G(X_\Sigma).$$

 This module is generated over $\hat
S(\g^*)$ by the elements $[\bar{O}_\sigma]$ for all $\sigma \in
\Sigma$. There are relations between these generators. However,
given  a scalar product on $\g$, we obtain  ${\rm
Todd}_{H_*^G}(X_\Sigma)$ as a combination of the classes
$[\bar{O}_\sigma]$, with canonical coefficients in $\hat S(\g^*)$ .

\begin{theorem}\label{todd}
Let $G$ be a torus with Lie algebra $\g$. Let $\Sigma$ be a
complete rational fan in $\g$ and let $X_\Sigma$ be the
corresponding $G$-toric variety. For $\sigma\in \Sigma$, let
$\bar{O}_\sigma$ be the corresponding $G$-invariant subvariety of
$X_\Sigma$, and let $[\bar{O}_\sigma]$ be its class in the
equivariant homology ring $H^G_*(X_\Sigma) $. Choose a scalar
product on $\g$. Then the equivariant Todd class ${\rm
Todd}_{H_*^G}(X_\Sigma)$ of $X_\Sigma$ is given by the following
combination of the elements $[\bar{O}_\sigma]$
\begin{equation}\label{toddexpansion}
{\rm Todd}_{H_*^G}(X_\Sigma)= \sum_{\sigma\in
\Sigma}\mu^*(\sigma)[\bar{O}_\sigma].
\end{equation}
\end{theorem}
\begin{proof}
We use the  push-forward argument of  \cite{brylinskizhang}.

When $\Sigma$ is a complete simplicial fan, the cohomological Todd
class is the Poincar{\'e} dual of the homological Todd class. Thus
Theorem \ref{todd}  follows right away from Proposition
\ref{expansionsimplicial} and Proposition \ref{FormulaP}.

For the general case, we  take advantage of the functoriality of the
Todd genus and use a refinement of the fan $\Sigma$, in order to
reduce the proof to the case of a complete simplicial toric variety.
The valuation property (Theorem \ref{valuation}) will be crucial.

The cones $\sigma \in \Sigma$  can be subdivided, yielding a fan
$\tilde{\Sigma}$, in such a way that all the cones in
$\tilde{\Sigma}$ are simplicial. Let $\tilde{X}$ be the toric
variety corresponding to the fan $\tilde{\Sigma}$. Then
$\tilde{X}$ is complete simplicial, and there is a proper
birational map $f: \tilde{X}\to X_\Sigma$. The equivariant  Todd
genus of $X_\Sigma$ coincides with the push-forward of the
equivariant Todd genus of $\tilde{X}$. Moreover, let
$\tilde{\sigma}$ be a cone in $\tilde{\Sigma}$, let
$\bar{O}_{\tilde{\sigma}}$ be the corresponding subvariety of
$\tilde{X}$, and let  $\sigma$ be the smallest cone of $\Sigma$
containing $\tilde{\sigma}$. Then
$f_*([\bar{O}_{\tilde{\sigma}}])= [\bar{O}_\sigma]$ if
$\tilde{\sigma}$ and $\sigma$ have the same dimension, and
$f_*([\bar{O}_{\tilde{\sigma}}])= 0$ otherwise.

For the  variety $\tilde{X}$, we have
\begin{equation}\label{smoothtoddexpansion}
{\rm Todd}_{H_*^G}(\tilde{X})= \sum_{\tilde{\sigma}\in
\tilde{\Sigma}}\mu^*(\tilde{\sigma})[\bar{O}_{\tilde{\sigma}}].
\end{equation}
 Let $\sigma\in\Sigma$. The
cones $\tilde{\sigma}\in\tilde{\Sigma}$ which are contained in
$\sigma$ form a subdivision of $\sigma$. Therefore we have by
Theorem \ref{valuation},
$$
\mu^*(\sigma)= \sum_{\{ \tilde{\sigma}\in\tilde{\Sigma} ,\;
\tilde{\sigma}\subset \sigma, \; \dim\tilde{\sigma}= \dim \sigma
\}}\mu^*(\tilde{\sigma}).
$$
Thus we obtain (\ref{toddexpansion}) by taking the push-forward of
both sides of (\ref{smoothtoddexpansion}). $\Box$ \end{proof}

\bigskip

By localizing Equation (\ref{toddexpansion}), we recover the formula
of Brylinski-Zhang (\cite{brylinskizhang}, Theorem 9.4), in the
following corollary.

We denote by $L$ the multiplicative subset of $S(\g^*)$ which
consists of products $\prod_{i\in I }v_i$ of nonzero elements of
$\g^*$. We denote by $L^{-1}H_*^G(X_\Sigma)$ the corresponding
localized module. It is generated by the classes of the fixed points
$p_{\sigma_0}$, for $\sigma_0\in\Sigma[d]$. More precisely we have
the following relations
\begin{lemma}\label{localize}
Let $\sigma$ be any cone of $\Sigma$. In $L^{-1}H_*^G(X_\Sigma)$, we
have
$$
[\bar{O}_\sigma]= \sum_{\{\sigma_0\in\Sigma[d], \sigma\in
\CF(\sigma_0)\}} I(\sigma_0^*\cap \lin(\sigma)^\perp)[p_{\sigma_0}].
$$
\end{lemma}
\begin{proof} The proof is by induction on the codimension of $\sigma$. If
$\sigma$ is of dimension $d$, there is nothing to prove. Otherwise,
take a nonzero $v\in \lin(\sigma)^\perp$. By Lemma 8.8 of
\cite{brylinskizhang}, we have
$$
-v [\bar{O}_\sigma]= \sum_{\{\tau\in\Sigma,
\sigma\in\CF'(\tau)\}}\langle \eta(\sigma,\tau),v \rangle
[\bar{O}_\tau]
$$
where $\CF'(\tau)$ denotes the set of facets of $\tau$ and
$\eta(\sigma,\tau)$ is the generator of the semigroup
$\tau\cap\Lambda/\sigma\cap\Lambda$.

Applying  the induction hypothesis to each $\tau$,  we have to
prove, for a fixed $\sigma_0\in \Sigma[d]$
\begin{equation}\label{stokes}
-v I(\sigma_0^*\cap \lin(\sigma)^\perp)=
\sum_{\{\tau\in\CF(\sigma_0), \sigma\in\CF'(\tau)\}}\langle
\eta(\sigma,\tau),v \rangle I(\sigma_0^*\cap \lin(\tau)^\perp).
\end{equation}

When $\tau$ runs through  the set of faces of $\sigma_0$ having
$\sigma$ as a facet, the cone $\sigma_0^*\cap \lin(\tau)^{\perp}$
runs through the set of facets of  $\sigma_0^*\cap
\lin(\sigma)^{\perp}$. Thus relation (\ref{stokes}) follows from
Lemma \ref{relationsstokes}. $\Box$
\end{proof}

\begin{corollary}
In $L^{-1}\hat{H}_*^G(X_\Sigma)$ the Todd class is given by
\begin{equation}\label{bryz}
 {\rm
Todd}_{H_*^G}(X_\Sigma)= \sum_{\sigma_0\in
\Sigma[d]}S(\sigma_0^*)[p_{\sigma_0}].
\end{equation}
\end{corollary}
\begin{proof}
We rewrite  (\ref{toddexpansion}) using Lemma \ref{localize} and
reverse the  summations.
 Let $\sigma_0\in \Sigma[d]$. The faces $\f$ of $\sigma_0^*$
 are in one to one correspondance with
 the faces $\sigma$  of $\sigma_0$, with $\f = \sigma_0^*\cap
 \lin(\sigma)^\perp$ corresponding to $\sigma$. Moreover, the
 transverse cone $\t(\sigma_0^*,\f)$ is equal to $\check{\sigma}$.
We thus obtain
$$
{\rm Todd}_{H_*^G}(X_\Sigma)= \sum_{\sigma_0\in \Sigma[d]} \sum
_{\f\in \CF(\sigma_0)} \mu(\t(\sigma_0^*,\f))I(\f)[p_{\sigma_0}].
$$
Hence the corollary follows from Euler-Maclaurin expansion
(\ref{definitionmu})  of $S(\sigma_0^*)$. $\Box$
\end{proof}

\vspace{2cm}
\thanks{Nicole Berline,
Ecole Polytechnique,  Centre de math\'ematiques Laurent Schwartz,
91128, Palaiseau, France.}

\thanks{email: berline@math.polytechnique.fr}

\vspace{1cm}
\thanks{Mich{\`e}le Vergne, Institut de Math\'ematiques de Jussieu, Th{\'e}orie des
Groupes, Case 7012, 2 Place Jussieu, 75251 Paris Cedex 05,
France;}

\thanks
{Ecole Polytechnique,  Centre de math\'ematiques Laurent Schwartz,
91128, Palaiseau, France.}

\thanks{email: vergne@math.polytechnique.fr}

\end{document}